\documentclass{amsart}

\newtheorem{theorem}{Theorem}[section]

\newtheorem{proposition}[theorem]{Proposition}
\newtheorem{corollary}[theorem]{Corollary}
\newtheorem{remark}[theorem]{Remark}
\theoremstyle{definition}
\newtheorem{definition}[theorem]{Definition}

\numberwithin{equation}{section}

\begin{document}
\title{Mabuchi and Aubin-Yau functionals over complex
surfaces}
\author{Yi Li}
\address{Department of Mathematics, Harvard University, Cambridge, MA 02138}

\email{yili@math.harvard.edu}

\dedicatory{This paper is dedicated to my colleague and friend Lin
Chen \\ who has passed away in an accident.}

\begin{abstract} In this note we construct Mabuchi
$\mathcal{L}^{{\rm M}}_{\omega}$ functional and Aubin-Yau
functionals $\mathcal{I}^{{\rm AY}}_{\omega}, \mathcal{J}^{{\rm
AY}}_{\omega}$ on any compact complex surfaces, and establish a
number of properties. Our construction coincides with the original
one in the K\"ahler case.
\end{abstract}
\maketitle

\tableofcontents
\section{Introduction}
Let $(X,\omega)$ be a compact K\"ahler manifold of the complex
dimension $n$. It's known that the volume $V_{\omega}$ depends only
on the K\"ahler class of $\omega$, namely,
\begin{equation*}
\int_{X}(\omega+\sqrt{-1}\partial\overline{\partial}\varphi)^{n}=\int_{X}\omega^{n}
\end{equation*}
for any real-valued smooth function $\varphi$ with
$\omega_{\varphi}:=\omega+\sqrt{-1}\partial\overline{\partial}\varphi>0$,
because of closedness of $\omega$.

If $(X,g)$ is a compact Hermitian manifold of the complex dimension
$n$, the same result does not hold in general. In the Section
\ref{Sec2}, we consider a function to describe such a phenomena,
i.e., we define
\begin{equation*}
{\rm
Err}_{\omega}(\varphi):=\int_{X}\omega^{n}-\int_{X}\omega^{n}_{\varphi},
\end{equation*}
where $\omega$ is its associated real $(1,1)$-form. The main result
\footnote{V. Tosatti told the author that the same result has
implicitly contained in \cite{TW2} and \cite{DK}.} is

\begin{proposition} \label{1.1} Let $(X,g)$ be a compact Hermitian manifold of the complex dimension $n$ and $\omega$ its associated
real $(1,1)$-form. If $\partial\overline{\partial}(\omega^{k})=0$
for $k=1,2$, then
\begin{equation*}
\int_{X}\omega^{n}_{\varphi}=\int_{X}\omega^{n}=:V_{\omega}
\end{equation*}
for any real-valued function $\varphi\in C^{\infty}(X)_{\mathbb{R}}$
with
$\omega_{\varphi}:=\omega+\sqrt{-1}\partial\overline{\partial}\varphi>0$.
\end{proposition}

In K\"ahler geometry, energy functionals, such as Mabuchi $K$-energy
functional \cite{M}, Aubin-Yau energy functionals \cite{PS}, and
Chen-Tian energy functionals \cite{CT}, play an important role in
studying K\"ahler-Einstein metrics and constant scalar curvatures.
When I was in Yau's Seminar, I asked myself that can we define
energy functionals on compact complex manifolds? This is one
motivation to write down this note. Another motivation comes from a
question in S.-T. Yau's survey \cite{Y2}, that find necessary and
sufficient conditions for a complex manifold to admit a K\"ahler
structure. When $n=2$, it was settled by Siu \cite{Siu} or see
\cite{LT}: A compact complex surface is K\"ahler if and only if its
first Betti number is even. In the second part of this note we
construct Mabuchi $\mathcal{L}^{{\rm M}}_{\omega}$ functional and
Aubin-Yau functionals $\mathcal{I}^{{\rm AY}}_{\omega},
\mathcal{J}^{{\rm AY}}_{\omega}$ on any compact complex surface.

Let $(X,g)$ be a compact Hermitian manifold of the complex dimension
$2$ and $\omega$ be its associated real $(1,1)$-form. Let
\begin{equation*}
\mathcal{P}_{\omega}:=\{\varphi\in C^{\infty}(X)_{\mathbb{R}}|
\omega_{\varphi}:=\omega+\sqrt{-1}\partial\overline{\partial}\varphi>0\}.
\end{equation*}
For any $\varphi',\varphi''\in\mathcal{P}_{\omega}$, we define
\begin{eqnarray*}
\mathcal{L}^{{\rm
M}}_{\omega}(\varphi',\varphi'')&:=&\frac{1}{V_{\omega}}\int^{1}_{0}\int_{X}\dot{\varphi}_{t}\cdot\omega^{2}_{\varphi_{t}}dt
\\
&-&\frac{1}{V_{\omega}}\int^{1}_{0}\int_{X}\sqrt{-1}\partial\omega\wedge(\overline{\partial}\dot{\varphi}_{t}\cdot\varphi_{t})dt
\\
&+&\frac{1}{V_{\omega}}\int^{1}_{0}\int_{X}\sqrt{-1}\overline{\partial}\omega\wedge
(\partial\dot{\varphi}_{t}\cdot\varphi_{t})dt
\end{eqnarray*}
where $\varphi_{t}$ is any smooth path in $\mathcal{P}_{\omega}$
from $\varphi'$ to $\varphi''$. Also we set
\begin{equation*}
\mathcal{L}^{{\rm M}}_{\omega}(\varphi):=\mathcal{L}^{{\rm
M}}_{\omega}(0,\varphi).
\end{equation*}

\begin{theorem} \label{1.2} The functional $\mathcal{L}^{{\rm
M}}_{\omega}(\varphi',\varphi'')$ is independent of the choice of
the smooth path $\{\varphi_{t}\}_{0\leq t\leq 1}$ and satisfies the
$1$-cocycle condition. In particular
\begin{eqnarray*}
\mathcal{L}^{{\rm
M}}_{\omega}(\varphi)&=&\frac{1}{3V_{\omega}}\int_{X}\varphi(\omega^{2}+\omega\wedge\omega_{\varphi}+\omega^{2}_{\varphi})
\\
&+&\frac{1}{2V_{\omega}}\int_{X}\varphi[-\sqrt{-1}\partial\omega\wedge\overline{\partial}\varphi+\sqrt{-1}\overline{\partial}\omega\wedge\partial\varphi].
\end{eqnarray*}
Moreover, for any $\varphi\in\mathcal{P}_{\omega}$ and any constant
$C\in\mathbb{R}$, we have
\begin{equation*}
\mathcal{L}^{{\rm
M}}_{\omega}(\varphi,\varphi+C)=C\cdot\left(1-\frac{{\rm
Err}_{\omega}(\varphi)}{V_{\omega}}\right);
\end{equation*}
for any $\varphi_{1},\varphi_{2}\in\mathcal{P}_{\omega}$ and any
constant $C\in\mathbb{R}$, we have
\begin{equation*}
\mathcal{L}^{{\rm
M}}_{\omega}(\varphi_{1},\varphi_{2}+C)=\mathcal{L}^{{\rm
M}}_{\omega}(\varphi_{1},\varphi_{2})+C\cdot\left(1-\frac{{\rm
Err}_{\omega}(\varphi)}{V_{\omega}}\right).
\end{equation*}
\end{theorem}

In Section \ref{Sec4} , we construct Aubin-Yau functionals on
compact complex surfaces. Let $(X,g)$ be a compact Hermitian
manifold of thee complex dimension $2$ and $\omega$ be its
associated real $(1,1)$-form. Set
\begin{eqnarray*}
\mathcal{A}_{\omega}(\varphi)&:=&\frac{1}{2V_{\omega}}\int_{X}\varphi\cdot-\sqrt{-1}\partial\omega\wedge\overline{\partial}\varphi,
\\
\mathcal{B}_{\omega}(\varphi)&:=&\frac{1}{2V_{\omega}}\int_{X}\varphi\cdot\sqrt{-1}\overline{\partial}\omega\wedge\partial\varphi.
\end{eqnarray*}
We define
\begin{eqnarray*}
\mathcal{I}^{{\rm
AY}}_{\omega}(\varphi)&:=&\frac{1}{V_{\omega}}\int_{X}\varphi(\omega^{2}-\omega^{2}_{\varphi})+2\mathcal{A}_{\omega}(\varphi)+2\mathcal{B}_{\omega}(\varphi),
\\
\mathcal{J}^{{\rm
AY}}_{\omega}(\varphi)&:=&\frac{1}{V_{\omega}}\int^{1}_{0}\int_{X}\varphi(\omega^{2}-\omega^{2}_{s\cdot\varphi})ds+\mathcal{A}_{\omega}(\varphi)+\mathcal{B}_{\omega}(\varphi).
\end{eqnarray*}

\begin{theorem} \label{1.3} For any compact Hermitian manifold $(X,g)$ of
the complex dimension $2$, we have
\begin{equation*}
\frac{1}{3}\mathcal{I}^{{\rm
AY}}_{\omega}(\varphi)\leq\mathcal{J}^{{\rm
AY}}_{\omega}(\varphi)\leq\frac{2}{3}\mathcal{I}^{{\rm
AY}}_{\omega}(\varphi)
\end{equation*}
for any $\varphi\in\mathcal{P}_{\omega}$, where $\omega$ is its
associated real $(1,1)$-form.
\end{theorem}

We hope this exposition will give some ideas to study Yau's problem.
The author are concerning the construction of those functionals on
higher dimensional compact complex manifolds.
\\

{\bf Acknowledgements.} The author would like to thank Valentino
Tosatti who read this note and pointed out a serious mistake in the
first version. Furthermore, I also thank Chen-Yu Chi and Ming-Tao
Chuan for useful discussion.

\section{${\rm Err}_{\omega}$ map on complex manifolds}\label{Sec2}
Let $(X,g)$ be a compact Hermitian manifold of the complex dimension
$n$ and write the associated real $(1,1)$-form as
\begin{equation*}
\omega=\sqrt{-1}\cdot g_{i\overline{j}}\cdot dz^{i}\wedge
d\overline{z^{j}}.
\end{equation*}
Let $\mathcal{P}_{\omega}$ be the space of all real-valued smooth
functions $\varphi\in C^{\infty}(X)_{\mathbb{R}}$ so that
$\omega+\sqrt{-1}\partial\overline{\partial}\varphi$ is positive
definite on $X$:
\begin{equation*}
\mathcal{P}_{\omega}:=\{\varphi\in C^{\infty}(X)_{\mathbb{R}}|
\omega+\sqrt{-1}\partial\overline{\partial}\varphi>0\}.
\end{equation*}
Also we set
\begin{equation*}
\mathcal{P}^{0}_{\omega}:=\left\{\varphi\in\mathcal{P}_{\omega}\Big|
\sup_{X}\varphi=0\right\}.
\end{equation*}
\subsection{${\rm Err}_{\omega}$ map on compact complex manifolds}
To such a function $\varphi\in\mathcal{P}_{\omega}$ we associate the
quantity
\begin{equation}
V_{\omega}(\varphi):=\int_{X}\omega^{n}_{\varphi} \label{2.1}
\end{equation}
the volume of $X$ with respect to $\varphi$. In particular we set
\begin{equation}
V_{\omega}:=V_{\omega}(0)=\int_{X}\omega^{n}.\label{2.2}
\end{equation}
For $X$ being K\"ahler, we have $V_{\omega}=V_{\omega}(\varphi)$ for
any $\varphi\in\mathcal{P}_{\omega}$. In the non-K\"ahler case, it's
not in general true. Hence it is reasonable to define
\begin{equation}
{\rm
Err}_{\omega}(\varphi):=V_{\omega}-V_{\omega}(\varphi)=\int_{X}\omega^{n}-\int_{X}\omega^{n}_{\varphi}.\label{2.3}
\end{equation}
A natural question is when does ${\rm Err}_{\omega}(\varphi)$ vanish
for any/a $\varphi\in\mathcal{P}_{\omega}$?Clearly there exists a
smooth real-valued function
$\varphi_{0}\equiv0\in\mathcal{P}_{\omega}$ such that
\begin{equation*}
\omega^{n}_{\varphi_{0}}=\omega^{n}, \ \ \ \sup_{X}\varphi_{0}=0,
\end{equation*}
hence
\begin{equation}
{\rm
Err}_{\omega}(\varphi_{0})=\int_{X}\omega^{n}-\int_{X}\omega^{n}=0.\label{2.4}
\end{equation}
This gives us some information about ${\rm Err}_{\omega}(\varphi)$
and motivates us to consider
\begin{equation}
{\rm
SupErr}_{\omega}:=\sup_{\varphi\in\mathcal{P}^{0}_{\omega}}\left({\rm
Err}_{\omega}(\varphi)\right), \ {\rm
InfErr}_{\omega}:=\inf_{\varphi\in\mathcal{P}^{0}_{\omega}}\left({\rm
Err}_{\omega}(\varphi)\right).\label{2.5}
\end{equation}
In any case, one has
\begin{equation}
{\rm InfErr}_{\omega}\leq0\leq{\rm
SupErr}_{\omega}\leq\int_{X}\omega^{n}.\label{2.6}
\end{equation}
It's interesting to find some conditions to guarantee that the
equalities hold. To study this behavior of ${\rm Err}_{\omega}$ we
consider the following several natural conditions on $\omega$:
\begin{itemize}

\item {\bf Condition 1.1:}
\begin{equation}
\sqrt{-1}\partial\omega\wedge\overline{\partial}\omega \ \
\text{and} \ \ \sqrt{-1}\partial\overline{\partial}\omega \ \
\text{are non-negative}\label{2.7}
\end{equation}

\item{\bf Condition 1.2:}

\begin{equation}
\sqrt{-1}\partial\omega\wedge\overline{\partial}\omega \ \
\text{and} \ \ \sqrt{-1}\partial\overline{\partial}\omega \ \
\text{are non-positive}\label{2.8}
\end{equation}

\item {\bf Condition 2:}
\begin{equation}
\partial\overline{\partial}(\omega^{k})=0, \ \ k=1,2.\label{2.9}
\end{equation}

\item {\bf Condition 3:}
\begin{equation}
d(\omega^{n-1})=0. \label{2.10}
\end{equation}

\item{\bf Condition 4:}
\begin{equation}
\partial\overline{\partial}(\omega^{n-1})=0.\label{2.11}
\end{equation}

\end{itemize}

\begin{remark} \label{2.1} Condition 2 was appeared in \cite{GL} as a
sufficient condition to solving the complex Monge-Amp\`ere equation
on Hermitian manifolds. The metric satisfying the third condition is
called a balanced metric, which naturally appears in string theory
(V. Tosatti and B. Wenkove \cite{TW1} solved the complex
Monge-Amp\`ere equation on Hermitian manifolds with balanced
metrics; later, they \cite{TW2} dropped off the balanced
condition.); When $n=2$, this condition is indeed the K\"ahler
condition. A metric satisfying Condition 4 is referred to be a
Gauduchon metric, and a theorem of Gauduchon \cite{G} shows that
there exists a Gauduchon metric on every compact Hermitian manifold.
Notice that Condition 3 implies condition 4, and Condition 2 is
equivalent to
$\partial\overline{\partial}\omega=0=\partial\omega\wedge\overline{\partial}\omega$.
In particular, Condition 2 implies Condition 1.1 and Condition 1.2.
In our case $n=2$, Condition 2 is equivalent to Condition 4.
\end{remark}

\begin{remark} \label{2.2} For any  two forms $\alpha$ of degree $|\alpha|$ and
$\beta$ of degree $|\beta|$, we have
\begin{equation*}
\alpha\wedge\beta=(-1)^{|\alpha|\cdot|\beta|}\beta\wedge\alpha.
\end{equation*}
Also we have
\begin{equation*}
d(\alpha\wedge\beta)=d\alpha\wedge\beta+(-1)^{|\alpha|}\alpha\wedge
d\beta;
\end{equation*}
according to types, one deduces
\begin{eqnarray*}
\partial(\alpha\wedge\beta)&=&\partial\alpha\wedge\beta+(-1)^{|\alpha|}\alpha\wedge\partial\beta,
\\
\overline{\partial}(\alpha\wedge\beta)&=&\overline{\partial}\alpha\wedge\beta+(-1)^{|\alpha|}\alpha\wedge\overline{\partial}\beta.
\end{eqnarray*}
Moreover, if $|\alpha|+|\beta|=2n-1$, then
\begin{eqnarray*}
\int_{X}\alpha\wedge\partial\beta&=&(-1)^{|\beta|}\int_{X}\partial\alpha\wedge\beta
\ \ = \ \ -(-1)^{|\alpha|}\int_{X}\partial\alpha\wedge\beta,
\\
\int_{X}\alpha\wedge\overline{\partial}\beta&=&(-1)^{|\beta|}\int_{X}\overline{\partial}\alpha\wedge\beta
\ \ = \ \
-(-1)^{|\alpha|}\int_{X}\overline{\partial}\alpha\wedge\beta.
\end{eqnarray*}
Those formulae are useful in our computations.

\end{remark}

\begin{theorem} \label{2.3} (i) If $\omega$ satisfies Condition 1.1, then
\begin{equation}
{\rm InfErr}_{\omega}=0.\label{2.12}
\end{equation}
(ii) Correspondingly, if $\omega$ satisfies Condition 1.2, then
\begin{equation}
{\rm SupErr}_{\omega}=0.\label{2.13}
\end{equation}
(iii) In particular ${\rm SupErr}_{\omega}={\rm InfErr}_{\omega}=0$
provided that $\omega$ satisfies Condition 2.
\end{theorem}

\begin{proof} (i) we knew that ${\rm
Err}_{\omega}(\varphi_{0})=0$ for some
$\varphi_{0}\in\mathcal{P}^{0}_{\omega}$. To achieve the argument,
we need only to show that ${\rm Err}_{\omega}(\varphi)\geq0$ for
each function $\varphi\in\mathcal{P}^{0}_{\omega}$. By definition we
have
\begin{eqnarray*}
& & {\rm
Err}_{\omega}(\varphi)\\
&=&-\int_{X}\omega^{n}_{\varphi}+\int_{X}\omega^{n} \ \ = \ \
\int_{X}-\sqrt{-1}\partial\overline{\partial}\varphi\wedge\sum^{n-1}_{i=0}\omega^{i}_{\varphi}\wedge\omega^{n-1-i}
\\
&=&\sum^{n-1}_{i=0}\int_{X}\omega^{i}_{\varphi}\wedge\omega^{n-1-i}\wedge(-\sqrt{-1}\partial\overline{\partial}\varphi)=\sum^{n-1}_{i=0}\int_{X}\sqrt{-1}\partial(\omega^{i}_{\varphi}\wedge\omega^{n-1-i})\wedge\overline{\partial}\varphi\\
&=&\sum^{n-1}_{i=0}\int_{X}\left[i\cdot\omega^{i-1}_{\varphi}
\wedge\partial\omega\wedge\omega^{n-1-i}+\omega^{i}_{\varphi}\wedge(n-1-i)\omega^{n-2-i}\wedge\partial\omega\right]\wedge\sqrt{-1}\overline{\partial}\varphi \\
&=&\sum^{n-1}_{i=0}\int_{X}\left[i\cdot\omega^{i-1}_{\varphi}\wedge\omega^{n-1-i}+(n-1-i)\omega^{i}_{\varphi}\wedge\omega^{n-2-i}\right]\wedge\partial\omega\wedge\sqrt{-1}\overline{\partial}\varphi
\\
&=&\sum^{n-1}_{i=0}\int_{X}\sqrt{-1}\left[\varphi\cdot(i\cdot\omega^{i-1}_{\varphi}\wedge\omega^{n-1-i}+(n-1-i)\cdot\omega^{i}_{\varphi}\wedge\omega^{n-2-i})\wedge\overline{\partial}\partial\omega
\right.\\
&+&\left.\varphi\cdot\overline{\partial}(i\cdot\omega^{i-1}_{\varphi}\wedge\omega^{n-1-i}+(n-1-i)\omega^{i}_{\varphi}\wedge\omega^{n-2-i})\wedge\partial\omega\right]
\\
&=&\sum^{n-1}_{i=0}(I_{i}+II_{i})
\end{eqnarray*}
where
\begin{eqnarray*}
I_{i}&=&\int_{X}\varphi\cdot[i\cdot\omega^{i-1}_{\varphi}\wedge\omega^{n-1-i}+(n-1-i)
\omega^{i}_{\varphi}\wedge\omega^{n-2-i}]\wedge(-\sqrt{-1}\partial\overline{\partial}\omega),
\\
II_{i}&=&\int_{X}\varphi\cdot\sqrt{-1}\overline{\partial}[i\cdot\omega^{i-1}_{\varphi}
\wedge\omega^{n-1-i}+(n-1-i)\omega^{i}_{\varphi}\wedge\omega^{n-2-i}]\wedge\partial\omega.
\end{eqnarray*}
Since $\sqrt{-1}\partial\overline{\partial}\omega\geq0$ and
$\varphi\leq0$ on $X$, the first term $I_{i}$ is non-negative.
Applying the integration by parts to $II_{i}$, we deduce
\begin{eqnarray*}
II_{i}&=&\int_{X}\varphi\left[i(i-1)\omega^{i-2}_{\varphi}\wedge\overline{\partial}\omega\wedge\omega^{n-1-i}\right.
+i(n-1-i)\omega^{i-1}_{\varphi}\wedge\omega^{n-2-i}\wedge\overline{\partial}\omega\\
\\
&+&\left.i(n-1-i)\omega^{i-1}_{\varphi}\wedge\overline{\partial}\omega\wedge\omega^{n-2-i}\right.
\\
&+&\left.(n-1-i)(n-2-i)\omega^{i}_{\varphi}\wedge\omega^{n-3-i}\wedge\overline{\partial}\omega\right]\wedge\sqrt{-1}\partial\omega
\\
&=&\int_{X}\varphi\cdot\omega^{i-2}_{\varphi}\wedge\omega^{n-3-i}\wedge[i(i-1)\omega^{2}+2i(n-1-i)\omega_{\varphi}\wedge\omega\\
&+&(n-1-i)(n-2-i)\omega^{2}_{\varphi}]
\wedge(-\sqrt{-1}\partial\omega\wedge\overline{\partial}\omega).
\end{eqnarray*}
Since $\sqrt{-1}\partial\omega\wedge\overline{\partial}\omega$ is
non-negative and $\varphi$ is non-positive, it follows that
$II_{i}\geq0$. Thus ${\rm Err}_{\omega}(\varphi)\geq0$ for each
$\varphi\in\mathcal{P}^{0}_{\omega}$ and therefore ${\rm
InfErr}_{\omega}=0$.

(ii) If $\omega$ satisfies Condition 1.2, the above proceeding gives
that ${\rm Err}_{\omega}(\varphi)\leq0$ for each
$\varphi\in\mathcal{P}^{0}_{\omega}$, i.e., ${\rm
SupErr}_{\omega}\leq0$. Hence ${\rm SupErr}_{\omega}=0$.

(iii) It's an immediate consequence of (i) and (ii).
\end{proof}

\begin{corollary} \label{2.4} If $\omega$ satisfies Condition 2, then ${\rm
Err}_{\omega}(\varphi)=0$ for any
$\varphi\in\mathcal{P}^{0}_{\omega}$. Equivalently, in this case,
the number $V_{\omega}(\varphi)=\int_{X}\omega^{n}_{\varphi}$ does
not depend on the choice of $\varphi\in\mathcal{P}_{\omega}$ and
equals $V_{\omega}=\int_{X}\omega^{n}$.
\end{corollary}

\subsection{Vanishing property of ${\rm Err}_{\omega}$ map on compact complex surface}
Let $(X,g)$ be a Hermitian manifold of the complex dimension $n$ and
let $\omega_{g}$ be its associated real $(1,1)$-form. We say that
$g$ is a Gauduchon metric if
$\partial\overline{\partial}(\omega^{n-1}_{g})=0$.

We recall a theorem of Gauduchon \cite{G} or see Remark \ref{2.1}.

\begin{theorem} \label{2.5} {\bf (Gauduchon, 1984)} If $X$ is a compact
complex manifold of the complex dimension $n$, then in the conformal
class of every Hermitian metric $g$ there exists a Gauduchon metric
$g_{{\rm G}}$, i.e., there is a positive function $\varphi\in
C^{\infty}(X)_{\mathbb{R}}$ such that $g_{{\rm G}}:=\varphi\cdot g$
is Gauduchon. If $X$ is connected and $n\geq2$, then $g_{{\rm G}}$
is unique up to a positive constants.
\end{theorem}

Using the existence of Gauduchon metric, we can prove the following

\begin{theorem} \label{2.6} Let $(X,g)$ be a compact complex surface with
Hermitian metric $g$ and let $\omega_{{\rm G}}$ be its associated
Gauduchon metric. Then
\begin{equation}
{\rm Err}_{\omega_{{\rm G}}}(\varphi)=0 \label{2.14}
\end{equation}
for all $\varphi\in\mathcal{P}_{\omega_{{\rm G}}}$.
\end{theorem}

\begin{proof} In what follows, we omit the subscript $G$ and write
$\omega_{{\rm G}}$ as $\omega$; since $\omega$ is a Gauduchon
metric, it follows that $\partial\overline{\partial}\omega=0$. In
the case of $n=2$, we have
\begin{eqnarray*}
{\rm
Err}_{\omega}(\varphi)&=&\int_{X}\omega^{2}-\int_{X}\omega^{2}_{\varphi}
\ \ = \ \
\int_{X}(\omega-\omega_{\varphi})(\omega+\omega_{\varphi}) \\
&=&\int_{X}(\omega+\omega_{\varphi})\wedge-\sqrt{-1}\partial\overline{\partial}\varphi
\\
&=&\int_{X}-\sqrt{-1}\partial\overline{\partial}(\omega+\omega_{\varphi})\cdot\varphi
\\
&=&\int_{X}-2\sqrt{-1}\partial\overline{\partial}\omega\cdot\varphi=0.
\end{eqnarray*}
Hence the theorem follows.
\end{proof}

\section{Mabuchi $\mathcal{L}^{{\rm M}}_{\omega}$ functional on complex surfaces}
\subsection{Mabuchi $\mathcal{L}^{{\rm M}}_{\omega}$ functional on compact
K\"ahler manifolds}\label{Sec3}
Suppose that $(X,\omega)$ is a compact K\"ahler manifold of the
complex dimension $n$. For any pair
$(\varphi',\varphi'')\in\mathcal{P}_{\omega}\times\mathcal{P}_{\omega}$
we define
\begin{equation*}
\mathcal{L}^{{\rm M}}_{\omega}:
\mathcal{P}_{\omega}\times\mathcal{P}_{\omega}\longrightarrow\mathbb{R}
\end{equation*}
as follows:
\begin{equation}
\mathcal{L}^{{\rm
M}}_{\omega}(\varphi',\varphi''):=\frac{1}{V_{\omega}}
\int^{1}_{0}\int_{X}\dot{\varphi}_{t}
\cdot\omega^{n}_{\varphi_{t}}\cdot dt\label{3.1}
\end{equation}
where $\{\varphi_{t}: 0\leq t\leq 1\}$ is any smooth path in
$\mathcal{P}_{\omega}$ such that $\varphi_{0}=\varphi'$ and
$\varphi_{1}=\varphi''$. For any $\varphi\in\mathcal{P}_{\omega}$ we
set
\begin{equation}
\mathcal{L}^{{\rm M}}_{\omega}(\varphi):=\mathcal{L}^{{\rm
M}}_{\omega}(0,\varphi).\label{3.2}
\end{equation}
Mabuchi \cite{M} showed that the functional (\ref{3.1}) is
well-defined, and hence we can explicitly write down
$\mathcal{L}^{{\rm M}}_{\omega}(\varphi)$.

In this section we extend Mabuchi $\mathcal{L}^{{\rm M}}_{\omega}$
functional to any compact complex surface by adding two extra terms
on the right hand side of (\ref{3.1}).
\subsection{Mabuchi $\mathcal{L}^{{\rm M}}_{\omega}$ functional on compact
complex surfaces}
Suppose now that $(X,g)$ is a compact complex surface and $\omega$
be its associated real $(1,1)$-form. Let
$\varphi',\varphi''\in\mathcal{P}_{\omega}$ and
$\{\varphi_{t}\}_{0\leq t\leq 1}$ be a smooth path in
$\mathcal{P}_{\omega}$ from $\varphi'$ to $\varphi''$.

Let
\begin{equation}
\mathcal{L}^{0}_{\omega}(\varphi',\varphi''):=\frac{1}{V_{\omega}}\int^{1}_{0}
\int_{X}\dot{\varphi}_{t}\cdot\omega^{2}_{\varphi_{t}}\cdot
dt.\label{3.3}
\end{equation}
Set
\begin{equation}
\psi(s,t):=s\cdot\varphi_{t}, \ \ \ 0\leq s\leq 1, \ \ \ 0\leq t\leq
1.\label{3.4}
\end{equation}
Consider a $1$-form on $[0,1]\times[0,1]$
\begin{equation}
\Psi^{0}:=\left(\int_{X}\frac{\partial\psi}{\partial
s}\cdot\omega^{2}_{\psi}\right)ds+\left(\int_{X}\frac{\partial\psi}{\partial
t}\cdot\omega^{2}_{\psi}\right)dt.\label{3.5}
\end{equation}
Taking differential on $\Psi^{0}$, we have
\begin{equation*}
d\Psi^{0}=I^{0}\cdot dt\wedge ds
\end{equation*}
where
\begin{equation}
I^{0}=\int_{X}\frac{\partial}{\partial
t}\left(\frac{\partial\psi}{\partial
s}\cdot\omega^{2}_{\psi}\right)-\int_{X}\frac{\partial}{\partial
s}\left(\frac{\partial\psi}{\partial
t}\cdot\omega^{2}_{\psi}\right).\label{3.6}
\end{equation}
Directly computing shows
\begin{eqnarray*}
I^{0}&=&\int_{X}\left[\frac{\partial^{2}\psi}{\partial t\partial
s}\cdot\omega^{2}_{\psi}+\frac{\partial\psi}{\partial
s}\cdot2\cdot\omega_{\psi}\wedge\sqrt{-1}\partial\overline{\partial}\left(\frac{\partial\psi}{\partial
t}\right)\right] \\
&-&\int_{X}\left[\frac{\partial^{2}\psi}{\partial s\partial
t}\cdot\omega^{2}_{\psi}+\frac{\partial\psi}{\partial t}\cdot 2\cdot
\omega_{\psi}\wedge\sqrt{-1}\partial\overline{\partial}\left(\frac{\partial\psi}{\partial
s}\right)\right] \\
&=&\int_{X}2\frac{\partial\psi}{\partial
s}\cdot\omega_{\psi}\wedge\sqrt{-1}\partial\overline{\partial}\left(\frac{\partial\psi}{\partial
t}\right)-\int_{X}2\frac{\partial\psi}{\partial
t}\cdot\omega_{\psi}\wedge\sqrt{-1}\partial\overline{\partial}\left(\frac{\partial\psi}{\partial
s}\right).
\end{eqnarray*}
In the following we deduce two slightly different formulae of
$I^{0}$. The first one is
\begin{eqnarray*}
I^{0}&=&\int_{X}2\frac{\partial\psi}{\partial
s}\cdot\omega_{\psi}\wedge\sqrt{-1}\partial\overline{\partial}\left(\frac{\partial\psi}{\partial
t}\right)+\int_{X}2\frac{\partial\psi}{\partial
t}\cdot\omega_{\psi}\wedge\sqrt{-1}\overline{\partial}\partial\left(\frac{\partial\psi}{\partial
s}\right) \\
&=&\int_{X}-2\sqrt{-1}\partial\left(\frac{\partial\psi}{\partial
s}\cdot\omega_{\psi}\right)\wedge\overline{\partial}\left(\frac{\partial\psi}{\partial
t}\right)\\
&+&\int_{x}-2\sqrt{-1}\overline{\partial}\left(\frac{\partial\psi}{\partial
t}\cdot\omega_{\psi}\right)\wedge\partial\left(\frac{\partial\psi}{\partial
s}\right) \\
&=&\int_{X}-2\sqrt{-1}\left[\partial\left(\frac{\partial\psi}{\partial
s}\right)\wedge\omega_{\psi}+\frac{\partial\psi}{\partial
s}\cdot\partial\omega\right]\wedge\overline{\partial}\left(\frac{\partial\psi}{\partial
t}\right) \\
&+&\int_{X}-2\sqrt{-1}\left[\overline{\partial}\left(\frac{\partial\psi}{\partial
t}\right)\wedge\omega_{\psi}+\frac{\partial\psi}{\partial
t}\cdot\overline{\partial}\omega\right]\wedge\partial\left(\frac{\partial\psi}{\partial
s}\right) \\
&=&\int_{X}-2\sqrt{-1}\frac{\partial\psi}{\partial
s}\cdot\partial\omega\wedge\overline{\partial}\left(\frac{\partial\psi}{\partial
t}\right)+\int_{X}-2\sqrt{-1}\frac{\partial\psi}{\partial
t}\cdot\overline{\partial}\omega\wedge\partial\left(\frac{\partial\psi}{\partial
s}\right) \\
&=&\int_{X}2\sqrt{-1}\frac{\partial\psi}{\partial
s}\cdot\overline{\partial}\left(\frac{\partial\psi}{\partial
t}\right)\wedge\partial\omega+\int_{X}2\sqrt{-1}\frac{\partial\psi}{\partial
t}\cdot\partial\left(\frac{\partial\psi}{\partial
s}\right)\wedge\overline{\partial}\omega.
\end{eqnarray*}
Similarly, we have
\begin{equation*}
I^{0}=\int_{X}-2\sqrt{-1}\frac{\partial\psi}{\partial
s}\cdot\partial\left(\frac{\partial\psi}{\partial
t}\right)\wedge\overline{\partial}\omega+\int_{X}-2\sqrt{-1}\frac{\partial\psi}{\partial
t}\cdot\overline{\partial}\left(\frac{\partial\psi}{\partial
s}\right)\wedge\partial\omega.
\end{equation*}

Next, we define
\begin{eqnarray}
\mathcal{L}^{1}_{\omega}(\varphi',\varphi'')&=&
\frac{1}{V_{\omega}}\int^{1}_{0}\int_{X}a_{2}\cdot\partial\omega\wedge(\overline{\partial}\dot{\varphi}_{t}\cdot\varphi_{t})dt,
\label{3.7}\\
\mathcal{L}^{2}_{\omega}(\varphi',\varphi'')&=&
\frac{1}{V_{\omega}}\int^{1}_{0}\int_{X}b_{2}\cdot\overline{\partial}\omega\wedge(\partial\dot{\varphi}_{t}\cdot\varphi_{t})dt.\label{3.8}
\end{eqnarray}
Here we require $\overline{a_{2}}=b_{2}$, and $a_{2}, b_{2}$ are
determined later. As before, consider
\begin{eqnarray*}
\Psi^{1}&=&\left[\int_{X}a_{2}\partial\omega\wedge\left(\overline{\partial}\left(\frac{\partial\psi}{\partial
s}\right)\cdot\psi\right)\right]ds+\left[\int_{X}a_{2}\partial\omega\wedge\left(\overline{\partial}\left(\frac{\partial\psi}{\partial
t}\right)\cdot\psi\right)\right]dt, \\
\Psi^{2}&=&\left[\int_{X}b_{2}\overline{\partial}\omega\wedge\left(\partial\left(\frac{\partial\psi}{\partial
s}\right)\cdot\psi\right)\right]ds+\left[\int_{X}b_{2}\overline{\partial}\omega\wedge\left(\partial\left(\frac{\partial\psi}{\partial
t}\right)\cdot\psi\right)\right]dt.
\end{eqnarray*}
Therefore
\begin{equation*}
d\Psi^{1}=I^{1}\cdot dt\wedge ds
\end{equation*}
where
\begin{eqnarray}
I^{1}&=&\int_{X}a_{2}\frac{\partial}{\partial
t}\left[\partial\omega\wedge\left(\overline{\partial}\left(\frac{\partial\psi}{\partial
s}\right)\cdot\psi\right)\right]
\label{3.9}\\
&-&\int_{X}a_{2}\frac{\partial}{\partial
s}\left[\partial\omega\wedge\left(\overline{\partial}\left(\frac{\partial\psi}{\partial
t}\right)\cdot\psi\right)\right].\nonumber
\end{eqnarray}
Dividing $I^{1}$ by $a_{2}$ yields
\begin{eqnarray*}
\frac{I^{1}}{a_{2}}&=&\int_{X}-\frac{\partial}{\partial
t}\left[\left(\overline{\partial}\left(\frac{\partial\psi}{\partial
s}\right)\cdot\psi\right)\wedge\partial\omega\right]+\int_{X}\frac{\partial}{\partial
s}\left[\left(\overline{\partial}\left(\frac{\partial\psi}{\partial
t}\right)\cdot\psi\right)\wedge\partial\omega\right] \\
&=&\int_{X}-\left[\overline{\partial}\left(\frac{\partial^{2}\psi}{\partial
t\partial
s}\right)\cdot\psi+\overline{\partial}\left(\frac{\partial\psi}{\partial
s}\right)\cdot\frac{\partial\psi}{\partial
t}\right]\wedge\partial\omega \\
&+&\int_{X}\left[\overline{\partial}\left(\frac{\partial^{2}\psi}{\partial
s\partial
t}\right)\cdot\psi+\overline{\partial}\left(\frac{\partial\psi}{\partial
t}\right)\cdot\frac{\partial\psi}{\partial
s}\right]\wedge\partial\omega \\
&=&\int_{X}-\frac{\partial\psi}{\partial
t}\cdot\overline{\partial}\left(\frac{\partial\psi}{\partial
s}\right)\wedge\partial\omega+\int_{X}\frac{\partial\psi}{\partial
s}\cdot\overline{\partial}\left(\frac{\partial\psi}{\partial
t}\right)\wedge\partial\omega.
\end{eqnarray*}
In the same way, one deduces
\begin{equation*}
d\Psi^{2}=I^{2}\cdot dt\wedge ds,
\end{equation*}
and
\begin{equation*}
\frac{I^{2}}{b_{2}}=\int_{X}-\frac{\partial\psi}{\partial
t}\cdot\partial\left(\frac{\partial\psi}{\partial
s}\right)\wedge\overline{\partial}\omega+\int_{X}\frac{\partial\psi}{\partial
s}\cdot\partial\left(\frac{\partial\psi}{\partial
t}\right)\wedge\overline{\partial}\omega.
\end{equation*}
Combining above formulas, we have
\begin{equation}
-\frac{I^{1}}{a_{2}}+\frac{I^{2}}{b_{2}}=-\frac{I^{0}}{\sqrt{-1}}.\label{3.10}
\end{equation}
Setting $a_{2}=-\sqrt{-1}$ and $b_{2}=\sqrt{-1}$, we get
\begin{equation*}
I^{0}+I^{1}+I^{2}=0.
\end{equation*}
Thus
\begin{equation}
d\Psi=0\label{3.11}
\end{equation}
where
\begin{equation*}
\Psi:=\Psi^{0}+\Psi^{1}+\Psi^{2}.
\end{equation*}

The following theorem is an immediate consequence of the above
discussion.

\begin{theorem} \label{3.1} Let $(X,g)$ be a compact complex surface and $\omega$ be its associated
real $(1,1)$-form. The functional
\begin{eqnarray}
\mathcal{L}^{{\rm
M}}_{\omega}(\varphi',\varphi'')&:=&\frac{1}{V_{\omega}}\int^{1}_{0}\int_{X}\dot{\varphi}_{t}\cdot\omega^{2}_{\varphi_{t}}dt
\label{3.12}\\
&-&\frac{1}{V_{\omega}}\int^{1}_{0}\int_{X}\sqrt{-1}\partial\omega\wedge(\overline{\partial}\dot{\varphi}_{t}\cdot\varphi_{t})dt
\nonumber\\
&+&\frac{1}{V_{\omega}}\int^{1}_{0}\int_{X}\sqrt{-1}\overline{\partial}\omega\wedge(\partial\dot{\varphi}_{t}\cdot\varphi_{t})dt\nonumber
\end{eqnarray}
is independent of the choice of the smooth path
$\{\varphi_{t}\}_{0\leq t\leq1}$. In particular,
\begin{eqnarray}
\mathcal{L}^{{\rm
M}}_{\omega}(\varphi)&:=&\mathcal{L}_{\omega}(0,\varphi) \ \ = \ \
\frac{1}{3V_{\omega}}\int_{X}\varphi(\omega^{2}+\omega\wedge\omega_{\varphi}+\omega^{2}_{\varphi})
\label{3.13}\\
&+&\frac{1}{2V_{\omega}}\int_{X}\varphi[-\sqrt{-1}\partial\omega\wedge\overline{\partial}\varphi+\sqrt{-1}\overline{\partial}\omega\wedge\partial\varphi].\nonumber
\end{eqnarray}
\end{theorem}

\begin{proof} Applying Stocks' theorem
to the region $\Delta=\{(s,t)\in\mathbb{R}^{2}: 0\leq s,t\leq 1\}$
and using equation (\ref{3.9}), we have
\begin{eqnarray*}
0&=&\int_{\Delta}d\Psi \ \ = \ \ \int_{\partial\Delta}\Psi \ \ = \ \
\int_{\partial\Delta}(\Psi^{0}+\Psi^{1}+\Psi^{2}) \\
&=&\int^{1}_{0}\dot{\varphi}_{t}\cdot\omega^{2}_{\varphi_{t}}dt-\int^{1}_{0}\int_{X}\frac{\partial\psi}{\partial
s}\cdot\omega^{2}_{\psi}ds\Big|^{t=1}_{t=0} \\
&+&\int^{1}_{0}\int_{X}-\sqrt{-1}\partial\omega\wedge(\overline{\partial}\dot{\varphi}_{t}\cdot\varphi_{t})dt
\\
&-&\int^{1}_{0}\int_{X}-\sqrt{-1}\partial\omega\wedge\left(\overline{\partial}\left(\frac{\partial\psi}{\partial
s}\right)\cdot\psi\right)ds\Big|^{t=1}_{t=0} \\
&+&\int^{1}_{0}\int_{X}\sqrt{-1}\overline{\partial}\omega\wedge(\partial\dot{\varphi}_{t}\cdot\varphi_{t})dt
\\
&-&\int^{1}_{0}\int_{X}\sqrt{-1}\overline{\partial}\omega\wedge\left(\partial\left(\frac{\partial\psi}{\partial
s}\right)\cdot\psi\right)ds\Big|^{t=1}_{t=0}.
\end{eqnarray*}
Equivalently,
\begin{eqnarray*}
\mathcal{L}^{{\rm
M}}_{\omega}(\varphi',\varphi'')&=&\int^{1}_{0}\int_{X}\frac{\partial\psi}{\partial
s}\cdot\omega^{2}_{\psi}ds\Big|^{t=1}_{t=0} \\
&+&\int^{1}_{0}\int_{X}-\sqrt{-1}\partial\omega\wedge\left(\overline{\partial}\left(\frac{\partial\psi}{\partial
s}\right)\cdot\psi\right)ds\Big|^{t=1}_{t=0} \\
&+&\int^{1}_{0}\int_{X}\sqrt{-1}\overline{\partial}\omega\wedge\left(\partial\left(\frac{\partial\psi}{\partial
s}\right)\cdot\psi\right)ds\Big|^{t=1}_{t=0}.
\end{eqnarray*}
It turns out that $\mathcal{L}^{{\rm
M}}_{\omega}(\varphi',\varphi'')$ is well-defined. For the second
argument, we can choose the smooth path $\varphi_{t}=t\cdot\varphi$,
$0\leq t\leq1$.
\end{proof}

\begin{remark}\label{3.2}  When $(X,g)$ is a compact K\"ahler surface, the
functional (\ref{3.12}) or (\ref{3.13}) coincides with the original
one.
\end{remark}

\begin{definition} \label{3.3} Suppose that $S$ is a non-empty set
and $A$ is an additive group. A mapping $\mathcal{N}: S\times S\to
A$ is said to satisfy the {\bf $1$-cocycle condition} if
\begin{itemize}

\item[(i)]
$\mathcal{N}(\sigma_{1},\sigma_{2})+\mathcal{N}(\sigma_{2},\sigma_{1})=0$;

\item[(ii)]
$\mathcal{N}(\sigma_{1},\sigma_{2})+\mathcal{N}(\sigma_{2},\sigma_{3})+\mathcal{N}(\sigma_{3},\sigma_{1})=0$.

\end{itemize}
\end{definition}

\begin{corollary} \label{3.4} The functional $\mathcal{L}^{{\rm M}}_{\omega}$ satisfies
the 1-cocycle condition.
\end{corollary}

\begin{corollary} \label{3.5} For any $\varphi\in\mathcal{P}_{\omega}$ and any
constant $C\in\mathbb{R}$, we have
\begin{equation}
\mathcal{L}^{{\rm
M}}_{\omega}(\varphi,\varphi+C)=C\cdot\left(1-\frac{{\rm
Err}_{\omega}(\varphi)}{V_{\omega}}\right). \label{3.14}
\end{equation}
In particular, if $\partial\overline{\partial}\omega=0$, then
$\mathcal{L}^{{\rm M}}_{\omega}(\varphi,\varphi+C)=C$.
\end{corollary}

\begin{proof} We choose the smooth path $\varphi_{t}=\varphi+t\cdot
C$, $t\in[0,1]$. So
\begin{eqnarray*}
\mathcal{L}^{{\rm
M}}_{\omega}(\varphi,\varphi+C)&=&\frac{1}{V_{\omega}}\int^{1}_{0}\int_{X}C\cdot\omega^{2}_{\varphi+C\cdot
t}dt \ \ = \ \
\frac{1}{V_{\omega}}\int^{1}_{0}\int_{X}C\cdot\omega^{2}_{\varphi}dt
\\
&=&\frac{1}{V_{\omega}}\int_{X}C\cdot\omega^{2}_{\varphi} \ \ = \ \
\frac{C}{V_{\omega}}\cdot V_{\omega}(\varphi) \ \ = \ \
C\cdot\left(1-\frac{{\rm Err}_{\omega}(\varphi)}{V_{\omega}}\right).
\end{eqnarray*}
If furthermore $\partial\overline{\partial}\omega=0$, then ${\rm
Err}_{\omega}(\varphi)=0$ for any $\varphi\in\mathcal{P}_{\omega}$.
\end{proof}

\begin{corollary} \label{3.6} For any
$\varphi_{1},\varphi_{2}\in\mathcal{P}_{\omega}$ and any constant
$C\in\mathbb{R}$, we have
\begin{equation}
\mathcal{L}^{{\rm
M}}_{\omega}(\varphi_{1},\varphi_{2}+C)=\mathcal{L}^{{\rm
M}}_{\omega}(\varphi_{1},\varphi_{2})+C\cdot\left(1-\frac{{\rm
Err}_{\omega}(\varphi_{2})}{V_{\omega}}\right).\label{3.2.13}
\end{equation}
\end{corollary}

\begin{proof} From Corollary \ref{3.4}, one has
\begin{equation*}
\mathcal{L}^{{\rm
M}}_{\omega}(\varphi_{1},\varphi_{2}+C)+\mathcal{L}^{{\rm
M}}_{\omega}(\varphi_{2},\varphi_{1}) =\mathcal{L}^{{\rm
M}}_{\omega}(\varphi_{2},\varphi_{2}+C).
\end{equation*}
Then the conclusion follows from Corollary 3.5.
\end{proof}
\section{Aubin-Yau functionals on compact complex
surfaces}\label{Sec4}
In this section we extend Aubin-Yau functionals to compact complex
surfaces, including K\"ahler surfaces, and deduce a number of basic
properties of these functionals.
\subsection{Aubin-Yau functionals on compact K\"ahler manifolds}
Suppose that $(X,\omega)$ is a compact K\"ahler manifold of
dimension $n$. For
$(\varphi',\varphi'')\in\mathcal{P}_{\omega}\times\mathcal{P}_{\omega}$,
Aubin-Yau functionals are defined by
\begin{eqnarray}
\mathcal{I}^{{\rm
AY}}_{\omega}(\varphi',\varphi'')&=&\frac{1}{V_{\omega}}
\int_{X}(\varphi''-\varphi')(\omega^{n}_{\varphi'}-\omega^{n}_{\varphi''}),\label{4.1}
\\
\mathcal{J}^{{\rm
AY}}_{\omega}(\varphi',\varphi'')&=&-\mathcal{L}^{{\rm
M}}_{\omega}(\varphi',\varphi'')+\frac{1}{V_{\omega}}\int_{X}
(\varphi''-\varphi')\omega^{n}_{\varphi'}.\label{4.2}
\end{eqnarray}
By definition, we have
\begin{equation}
\mathcal{J}^{{\rm
AY}}_{\omega}(\varphi',\varphi'')+\mathcal{J}^{{\rm
AY}}_{\omega}(\varphi'',\varphi')=\mathcal{I}^{{\rm
AY}}_{\omega}(\varphi',\varphi'')=\mathcal{I}^{{\rm
AY}}_{\omega}(\varphi'',\varphi').\label{4.3}
\end{equation}
For any $\varphi\in\mathcal{P}_{\omega}$ we set
\begin{eqnarray}
\mathcal{I}^{{\rm
AY}}_{\omega}(\varphi)&=&\frac{1}{V_{\omega}}\int_{X}\varphi(\omega^{n}-\omega^{n}_{\varphi}),
\label{4.4}\\
\mathcal{J}^{{\rm
AY}}_{\omega}(\varphi)&=&\int^{1}_{0}\frac{\mathcal{I}_{\omega}(s\cdot\varphi)}{s}ds=\frac{1}{V_{\omega}}\int^{1}_{0}\int_{X}\varphi(\omega^{n}-\omega^{n}_{s\cdot\varphi})ds.\label{4.5}
\end{eqnarray}
It's clear that
\begin{equation}
\mathcal{I}_{\omega}^{{\rm AY}}(\varphi)=\mathcal{I}^{{\rm
AY}}_{\omega}(0,\varphi), \ \ \ \mathcal{J}^{{\rm
AY}}_{\omega}(\varphi)=\mathcal{J}^{{\rm
AY}}_{\omega}(0,\varphi).\label{4.6}
\end{equation}
By definition we have
\begin{eqnarray}
& & \mathcal{J}^{{\rm AY}}_{\omega}(\varphi)\nonumber
\\
&=&\frac{1}{V_{\omega}}\int^{1}_{0}ds\int_{X}\varphi(-\sqrt{-1}\partial\overline{\partial}(s\cdot\varphi))\wedge\sum^{n-1}_{i=1}\omega^{n-1-i}\wedge\omega^{i}_{s\cdot\varphi}
\nonumber\\
&=&\frac{-\sqrt{-1}}{V_{\omega}}\int^{1}_{0}s\cdot
ds\int_{X}\varphi\cdot\partial\overline{\partial}\varphi\wedge\sum^{n-1}_{i=0}\omega^{n-1-i}\wedge[\omega+s(\omega_{\varphi}-\omega)]^{i}
\nonumber\\
&=&\frac{-\sqrt{-1}}{V_{\omega}}\int^{1}_{0}s\cdot
ds\int_{X}\varphi\cdot\partial\overline{\partial}\varphi\sum^{n-1}_{i=0}\omega^{n-1-i}\wedge\sum^{i}_{j=0}\binom{i}{j}(1-s)^{i-j}s^{j}\omega^{i-j}\wedge\omega^{j}_{\varphi}
\nonumber\\
&=&\frac{-\sqrt{-1}}{V_{\omega}}\int_{X}\varphi\cdot\partial\overline{\partial}\varphi\wedge\sum^{n-1}_{j=0}\sum^{n-1}_{i=j}\binom{i}{j}\omega^{n-1-j}\wedge\omega^{j}_{\varphi}\int^{1}_{0}(1-s)^{i-j}s^{1+j}
\nonumber\\
&=&\frac{-\sqrt{-1}}{V_{\omega}}\int_{X}\varphi\cdot\partial\overline{\partial}\varphi\wedge\sum^{n-1}_{j=0}\omega^{n-1-j}\wedge\omega^{j}_{\varphi}\cdot\sum^{n-1}_{i=j}\binom{i}{j}\cdot\frac{(i-j)!\cdot
(j+1)!}{(i+2)!} \nonumber\\
&=&\frac{-\sqrt{-1}}{V_{\omega}}\int_{X}\varphi\cdot\partial\overline{\partial}\varphi\wedge\sum^{n-1}_{j=0}\omega^{n-1-j}\wedge\omega^{j}_{\varphi}\cdot\sum^{n-1}_{i=j}\frac{j+1}{(i+2)(i+1)}
\nonumber\\
&=&\frac{-\sqrt{-1}}{V_{\omega}}\int_{X}\varphi\cdot\partial\overline{\partial}\varphi\wedge\sum^{n-1}_{j=0}\frac{n-j}{n+1}\omega^{n-1-j}\wedge\omega^{j}_{\varphi}\nonumber
\end{eqnarray}
since
\begin{equation*}
\sum^{n-1}_{i=j}\frac{1}{(i+2)(i+1)}=\sum^{n-1}_{i=j}\left(\frac{1}{i+1}-\frac{1}{i+2}\right)=\frac{1}{j+1}-\frac{1}{n+1}.
\end{equation*}
On the other hand,
\begin{equation*}
\frac{n}{n+1}\mathcal{I}^{{\rm
AY}}_{\omega}(\varphi)=\frac{-\sqrt{-1}}{V_{\omega}}\int_{X}\varphi\cdot\partial\overline{\partial}\varphi\wedge\sum^{n-1}_{i=0}\frac{n}{n+1}\omega^{n-1-i}\wedge\omega^{i}_{\varphi}
\end{equation*}
Hence
\begin{eqnarray}
& & \frac{n}{n+1}\mathcal{I}^{{\rm
AY}}_{\omega}(\varphi)-\mathcal{J}^{{\rm AY}}_{\omega}(\varphi)
\label{4.7}\\
&=&\frac{1}{V_{\omega}}\int_{X}\varphi\cdot(-\sqrt{-1}\partial\overline{\partial}\varphi)\wedge\sum^{n-1}_{j=1}\frac{j}{n+1}\omega^{n-1-j}\wedge\omega^{j}_{\varphi}.\nonumber
\end{eqnarray}
Moreover,
\begin{eqnarray}
& & (n+1)\mathcal{J}^{{\rm AY}}_{\omega}(\varphi)-\mathcal{I}^{{\rm
AY}}_{\omega}(\varphi) \label{4.8}\\
&=&\frac{1}{V_{\omega}}\int_{X}\varphi\cdot-\sqrt{-1}\partial\overline{\partial}\varphi\wedge\sum^{n-1}_{j=0}(n-1-j)\omega^{n-1-j}\wedge\omega^{j}_{\varphi}.\nonumber
\end{eqnarray}

\begin{remark} \label{4.1} Notice that formulae (\ref{4.7}) and (\ref{4.8}) are also valid when $\omega$ is
non-K\"ahler.
\end{remark}

\subsection{Aubin-Yau functionals over compact complex surfaces}
Let $(X,g)$ be a compact complex manifold of the complex dimension
$n$ and $\omega$ be its associated real $(1,1)$-form. From Remark
4.1, we can formally use the notion $\mathcal{I}^{{\rm
AY}}_{\omega}$, $\mathcal{J}^{{\rm AY}}_{\omega}$, and
$\mathcal{L}^{{\rm M}}_{\omega}$, but now $\omega$ may not be
K\"ahler. Precisely, for any $\varphi\in\mathcal{P}_{\omega}$ we set
\begin{eqnarray}
\mathcal{I}^{{\rm
AY}}_{\omega|\bullet}(\varphi)&=&\frac{1}{V_{\omega}}\int_{X}\varphi(\omega^{n}-\omega^{n}_{\varphi}),
\label{4.9}\\
\mathcal{J}^{{\rm
AY}}_{\omega|\bullet}(\varphi)&=&\int^{1}_{0}\frac{\mathcal{I}^{{\rm
AY}}_{\omega|\bullet}(s\cdot\varphi)}{s}ds=\frac{1}{V_{\omega}}\int^{1}_{0}\int_{X}\varphi(\omega^{n}-\omega^{n}_{s\cdot\varphi})ds.\label{4.10}
\end{eqnarray}
Hence
\begin{eqnarray}
& & \frac{n}{n+1}\mathcal{I}^{{\rm
AY}}_{\omega|\bullet}(\varphi)-\mathcal{J}^{{\rm
AY}}_{\omega|\bullet}(\varphi)
\label{4.11}\\
&=&\frac{1}{V_{\omega}}\int_{X}\varphi\cdot(-\sqrt{-1}\partial\overline{\partial}\varphi)\wedge\sum^{n-1}_{j=1}\frac{j}{n+1}\omega^{n-1-j}\wedge\omega^{j}_{\varphi}.\nonumber
\end{eqnarray}
Moreover,
\begin{eqnarray}
& & (n+1)\mathcal{J}^{{\rm
AY}}_{\omega|\bullet}(\varphi)-\mathcal{I}^{{\rm
AY}}_{\omega|\bullet}(\varphi) \label{4.12}\\
&=&\frac{1}{V_{\omega}}\int_{X}\varphi\cdot-\sqrt{-1}\partial\overline{\partial}\varphi\wedge\sum^{n-1}_{j=0}(n-1-j)\omega^{n-1-j}\wedge\omega^{j}_{\varphi}.\nonumber
\end{eqnarray}

Restricting to compact complex surfaces and introducing two extra
functionals on $\mathcal{P}_{\omega}$
\begin{eqnarray}
\mathcal{A}_{\omega}(\varphi)&:=&\frac{1}{2V_{\omega}}\int_{X}\varphi\cdot-\sqrt{-1}\partial\omega\wedge\overline{\partial}\varphi,
\label{4.13}\\
\mathcal{B}_{\omega}(\varphi)&:=&\frac{1}{2V_{\omega}}\int_{X}\varphi\cdot\sqrt{-1}\overline{\partial}\omega\wedge\partial\varphi,\label{4.14}
\end{eqnarray}
(clearly
$\overline{\mathcal{A}_{\omega}(\varphi)}=\mathcal{B}_{\omega}(\varphi)$),
we define Aubin-Yau functionals as follows (Here constants $a,b,c,d$
are determined later, and actually $a=b=c=d=2$)
\begin{eqnarray}
\mathcal{I}^{{\rm AY}}_{\omega}(\varphi)&:=&\mathcal{I}^{{\rm
AY}}_{\omega|\bullet}(\varphi)+a\mathcal{A}_{\omega}(\varphi)+b\mathcal{B}_{\omega}(\varphi),
\label{4.15}\\
\mathcal{J}^{{\rm AY}}_{\omega}(\varphi)&:=&-\mathcal{L}^{{\rm
M}}_{\omega}(\varphi)+\frac{1}{V_{\omega}}\int_{X}\varphi\cdot\omega^{2}+c\mathcal{A}_{\omega}(\varphi)+d\mathcal{B}_{\omega}(\varphi).\label{4.16}
\end{eqnarray}
Since
\begin{eqnarray*}
\mathcal{J}^{{\rm
AY}}_{\omega|\bullet}(\varphi)&=&\frac{1}{V_{\omega}}\int_{X}\varphi\cdot(-\sqrt{-1}\partial\overline{\partial}\varphi)\wedge\sum^{1}_{j=0}\frac{2-j}{3}(\omega^{1-j}\wedge\omega^{j}_{\varphi})
\\
&=&\frac{1}{V_{\omega}}\int_{X}\varphi(\omega-\omega_{\varphi})\wedge\left(\frac{2}{3}\omega+\frac{1}{3}\omega_{\varphi}\right)
\\
&=&\frac{1}{V_{\omega}}\int_{X}\varphi\left(\frac{2}{3}\omega^{2}-\frac{1}{3}\omega\wedge\omega_{\varphi}-\frac{1}{3}\omega^{2}_{\varphi}\right),
\end{eqnarray*}
it follows that
\begin{equation}
\mathcal{J}^{{\rm AY}}_{\omega}(\varphi)=\mathcal{J}^{{\rm
AY}}_{\omega|\bullet}(\varphi)+(c-1)\mathcal{A}_{\omega}(\varphi)+(d-1)\mathcal{B}_{\omega}(\varphi).\label{4.17}
\end{equation}
Using Remark \ref{4.1} or the previous subsection, we deduce
\begin{eqnarray*}
& & \frac{2}{3}(\mathcal{I}^{{\rm
AY}}_{\omega}(\varphi)-a\mathcal{A}_{\omega}(\varphi)-b\mathcal{B}_{\omega}(\varphi))-\left(\mathcal{J}^{{\rm
AY}}_{\omega}(\varphi)-(c-1)\mathcal{A}_{\omega}(\varphi)-(d-1)\mathcal{B}_{\omega}(\varphi)\right)
\\
&=&\frac{1}{V_{\omega}}\int_{X}\varphi(-\sqrt{-1}\partial\overline{\partial}\varphi)\sum^{1}_{j=1}\frac{j}{3}\omega^{1-j}\wedge\omega^{j}_{\varphi}
\\
&=&\frac{1}{V_{\omega}}\int_{X}\varphi\cdot\frac{1}{3}\omega_{\varphi}\wedge(-\sqrt{-1}\partial\overline{\partial}\varphi)
\ \ = \ \
\frac{1}{V_{\omega}}\int_{X}\varphi\cdot\frac{1}{3}\omega_{\varphi}\wedge\sqrt{-1}\overline{\partial}\partial\varphi.
\end{eqnarray*}
Thus the left hand side has two slightly different expressions. If
we adopt the first one, we have
\begin{eqnarray*}
& & \frac{2}{3}(\mathcal{I}^{{\rm
AY}}_{\omega}(\varphi)-a\mathcal{A}_{\omega}(\varphi)-b\mathcal{B}_{\omega}(\varphi))-\left(\mathcal{J}^{{\rm
AY}}_{\omega}(\varphi)-(c-1)\mathcal{A}_{\omega}(\varphi)-(d-1)\mathcal{B}_{\omega}(\varphi)\right)
\\
&=&\frac{1}{3V_{\omega}}\int_{X}\sqrt{-1}\partial(\varphi\cdot\omega_{\varphi})\wedge\overline{\partial}\varphi
\ \ = \ \
\frac{1}{3V_{\omega}}\int_{X}\sqrt{-1}(\partial\varphi\wedge\omega_{\varphi}+\varphi\cdot\partial\omega)\wedge\overline{\partial}\varphi
\\
&=&\frac{1}{3V_{\omega}}\int_{X}\sqrt{-1}\partial\varphi\wedge\overline{\partial}\varphi\wedge\omega_{\varphi}-\frac{2}{3}\mathcal{A}_{\omega}(\varphi).
\end{eqnarray*}
On the other hand, using the second expression gives
\begin{eqnarray*}
& & \frac{2}{3}(\mathcal{I}^{{\rm
AY}}_{\omega}(\varphi)-a\mathcal{A}_{\omega}(\varphi)-b\mathcal{B}_{\omega}(\varphi))-\left(\mathcal{J}^{{\rm
AY}}_{\omega}(\varphi)-(c-1)\mathcal{A}_{\omega}(\varphi)-(d-1)\mathcal{B}_{\omega}(\varphi)\right)
\\
&=&\frac{1}{3V_{\omega}}\int_{X}-\sqrt{-1}\overline{\partial}(\varphi\cdot\omega_{\varphi})\wedge\partial\varphi
\ \ = \ \
\frac{1}{3V_{\omega}}\int_{X}-\sqrt{-1}(\overline{\partial}\varphi\wedge\omega_{\varphi}+\varphi\cdot\overline{\partial}\omega)\wedge\partial\varphi
\\
&=&\frac{1}{3V_{\omega}}\int_{X}\sqrt{-1}\partial\varphi\wedge\overline{\partial}\varphi\wedge\omega_{\varphi}-\frac{2}{3}\mathcal{B}_{\omega}(\varphi).
\end{eqnarray*}
Therefore
\begin{eqnarray*}
& & \frac{2}{3}(\mathcal{I}^{{\rm
AY}}_{\omega}(\varphi)-a\mathcal{A}_{\omega}(\varphi)-b\mathcal{B}_{\omega}(\varphi))-\left(\mathcal{J}^{{\rm
AY}}_{\omega}(\varphi)-(c-1)\mathcal{A}_{\omega}(\varphi)-(d-1)\mathcal{B}_{\omega}(\varphi)\right)
\\
&=&\frac{1}{3V_{\omega}}\int_{X}\sqrt{-1}\partial\varphi\wedge\overline{\partial}\varphi\wedge\omega_{\varphi}-\frac{\mathcal{A}_{\omega}(\varphi)+\mathcal{B}_{\omega}(\varphi)}{3},
\end{eqnarray*}
or, equivalently,
\begin{equation}
\frac{2}{3}\mathcal{I}^{{\rm
AY}}_{\omega}(\varphi)-\mathcal{J}^{{\rm
AY}}_{\omega}(\varphi)=\frac{1}{3V_{\omega}}\int_{X}\sqrt{-1}\partial\varphi\wedge\overline{\partial}\varphi\wedge\omega_{\varphi}\label{4.18}
\end{equation}
where we require
\begin{equation}
\frac{2}{3}a-(c-1)-\frac{1}{3}=0=\frac{2}{3}b-(d-1)-\frac{1}{3}.\label{4.19}
\end{equation}

\begin{theorem} \label{4.2} For any $\varphi\in\mathcal{P}_{\omega}$, one has
\begin{equation}
\frac{2}{3}\mathcal{I}^{{\rm
AY}}_{\omega}(\varphi)-\mathcal{J}^{{\rm
AY}}_{\omega}(\varphi)\geq0.\label{4.20}
\end{equation}
\end{theorem}

Using (\ref{4.12}) yields
\begin{eqnarray*}
& & 3\left(\mathcal{J}^{{\rm
AY}}_{\omega}(\varphi)-\frac{1}{2}\mathcal{A}_{\omega}(\varphi)-\frac{1}{2}\mathcal{B}_{\omega}(\varphi)\right)-(\mathcal{I}^{{\rm
AY}}_{\omega}(\varphi)-\mathcal{A}_{\omega}(\varphi)-\mathcal{B}_{\omega}(\varphi))
\\
&=&\frac{1}{V_{\omega}}\int_{X}\varphi\cdot(-\sqrt{-1}\partial\overline{\partial}\varphi)\wedge\sum^{1}_{j=0}(1-j)\omega^{1-j}\wedge\omega^{j}_{\varphi}
\\
&=&\frac{1}{V_{\omega}}\int_{X}\varphi\cdot(-\sqrt{-1}\partial\overline{\partial}\varphi)\wedge\omega
\\
&=&\frac{1}{V_{\omega}}\int_{X}(\varphi\cdot\omega)\wedge(-\sqrt{-1}\partial\overline{\partial}\varphi)
\ \ = \ \
\frac{1}{V_{\omega}}\int_{X}(\varphi\cdot\omega)\wedge\sqrt{-1}\overline{\partial}\partial\varphi.
\end{eqnarray*}
As the proof of Theorem \ref{4.2}, we have
\begin{eqnarray*}
& & 3\left(\mathcal{J}^{{\rm
AY}}_{\omega}(\varphi)-(c-1)\mathcal{A}_{\omega}(\varphi)-(d-1)\mathcal{B}_{\omega}(\varphi)\right)-(\mathcal{I}^{{\rm
AY}}_{\omega}(\varphi)-a\mathcal{A}_{\omega}(\varphi)-b\mathcal{B}_{\omega}(\varphi))
\\
&=&\frac{1}{V_{\omega}}\int_{X}\sqrt{-1}\partial(\varphi\cdot\omega)\wedge\overline{\partial}\varphi
\ \ = \ \
\frac{1}{V_{\omega}}\int_{X}\sqrt{-1}(\partial\varphi\wedge\omega+\varphi\cdot\partial\omega)\wedge\overline{\partial}\varphi
\\
&=&\frac{1}{V_{\omega}}\int_{X}\sqrt{-1}\partial\varphi\wedge\overline{\partial}\varphi\wedge\omega-2\mathcal{A}_{\omega}(\varphi)
\end{eqnarray*}
and
\begin{eqnarray*}
& & 3\left(\mathcal{J}^{{\rm
AY}}_{\omega}(\varphi)-(c-1)\mathcal{A}_{\omega}(\varphi)-(d-1)\mathcal{B}_{\omega}(\varphi)\right)-(\mathcal{I}^{{\rm
AY}}_{\omega}(\varphi)-a\mathcal{A}_{\omega}(\varphi)-b\mathcal{B}_{\omega}(\varphi))
\\
&=&\frac{1}{V_{\omega}}\int_{X}-\sqrt{-1}\overline{\partial}(\varphi\cdot\omega)\wedge\partial\varphi
\ \ = \ \
\frac{1}{V_{\omega}}\int_{X}-\sqrt{-1}(\overline{\partial}\varphi\wedge\omega+\varphi\cdot\overline{\partial}\omega)\wedge\partial\varphi
\\
&=&\frac{1}{V_{\omega}}\int_{X}\sqrt{-1}\partial\varphi\wedge\overline{\partial}\varphi\wedge\omega-2\mathcal{B}_{\omega}(\varphi).
\end{eqnarray*}
Hence
\begin{eqnarray*}
& & 3\left(\mathcal{J}^{{\rm
AY}}_{\omega}(\varphi)-(c-1)\mathcal{A}_{\omega}(\varphi)-(d-1)\mathcal{B}_{\omega}(\varphi)\right)-(\mathcal{I}^{{\rm
AY}}_{\omega}(\varphi)-a\mathcal{A}_{\omega}(\varphi)-b\mathcal{B}_{\omega}(\varphi))
\\
&=&\frac{1}{V_{\omega}}\int_{X}\sqrt{-1}\partial\varphi\wedge\overline{\partial}\varphi\wedge\omega
-(\mathcal{A}_{\omega}(\varphi)+\mathcal{B}_{\omega}(\varphi)).
\end{eqnarray*}
Equivalently,
\begin{equation}
3\mathcal{J}^{{\rm AY}}_{\omega}(\varphi)-\mathcal{I}^{{\rm
AY}}_{\omega}(\varphi)=\frac{1}{V_{\omega}}\int_{X}\sqrt{-1}\partial\varphi\wedge\overline{\partial}\varphi\wedge\omega,\label{4.21}
\end{equation}
where we also require
\begin{equation}
3(c-1)-a-1=0=3(d-1)-b-1.\label{4.22}
\end{equation}

\begin{theorem} \label{4.3} For any $\varphi\in\mathcal{P}_{\omega}$, one has
\begin{equation}
3\mathcal{J}^{{\rm AY}}_{\omega}(\varphi)-\mathcal{I}^{{\rm
AY}}_{\omega}(\varphi)\geq0. \label{4.23}
\end{equation}
\end{theorem}

Combining (\ref{4.19}) and (\ref{4.22}) we obtain the value of those
constants:
\begin{equation*}
a=b=c=d=2.
\end{equation*}

\begin{corollary} \label{4.4} For any compact complex surface $(X,g)$ and
any real-valued smooth function $\varphi\in\mathcal{P}_{\omega}$, we
have
\begin{eqnarray*}
& & \frac{1}{3}\mathcal{I}^{{\rm
AY}}_{\omega}(\varphi)\leq\mathcal{J}^{{\rm
AY}}_{\omega}(\varphi)\leq\frac{2}{3}\mathcal{I}^{{\rm
AY}}_{\omega}(\varphi), \\
& & \frac{3}{2}\mathcal{J}^{{\rm
AY}}_{\omega}(\varphi)\leq\mathcal{I}_{\omega}(\varphi)\leq
3\mathcal{J}^{{\rm AY}}_{\omega}(\varphi), \\
& & \frac{1}{2}\mathcal{J}^{{\rm
AY}}_{\omega}(\varphi)\leq\frac{1}{3}\mathcal{I}^{{\rm
AY}}_{\omega}(\varphi)\leq\mathcal{I}^{{\rm
AY}}_{\omega}(\varphi)-\mathcal{J}^{{\rm
AY}}_{\omega}(\varphi)\leq\frac{2}{3}\mathcal{I}^{{\rm
AY}}_{\omega}(\varphi)\leq 2\mathcal{J}^{{\rm
AY}}_{\omega}(\varphi),
\end{eqnarray*}
where $\omega$ is its associated real $(1,1)$-form.
\end{corollary}

\bibliographystyle{amsplain}

\end{document}